\documentclass[12pt,twoside]{amsart}
\usepackage{amsfonts}
\usepackage{amssymb}
\usepackage{amsmath}
\usepackage{latexsym}
\usepackage{enumerate}
\begin{document}

\title{Interpolation and Balls in $\mathbb{C}^k$}

\author[V.~I.~Paulsen]{Vern I.~Paulsen*}
\address{Dept.~of Mathematics, University of Houston, Houston, 
Texas
77204-3476,
U.S.A.}
\email{vern@math.uh.edu}
\thanks{*Research supported in part by a grant from the NSF}
\author[J.~P.~Solazzo]{James P.~Solazzo}
\address{Dept.~of Mathematics, Coastal Carolina University, Conway, SC 
29526-6054,U.S.A.}
\email{jsolazzo@coastal.edu}

\maketitle
\pagestyle{headings}
\markboth{Vern I. Paulsen and James P. Solazzo }{Interpolation 
and
Balls in $\mathbb{C}^k$}

\baselineskip=20pt

\theoremstyle{plain}
\newtheorem{theorem}{Theorem}[section]
\newtheorem{proposition}[theorem]{Proposition}
\newtheorem{lemma}[theorem]{Lemma}
\newtheorem{corollary}[theorem]{Corollary}
\newtheorem{problem}[theorem]{Problem}
\newtheorem{conjecture}[theorem]{Conjecture}
\newtheorem{definition}[theorem]{Definition}
\newtheorem*{remark}{Remark}
\newtheorem{example}[theorem]{Example}
\newcommand{\pa}{\ensuremath{{\mathcal{P}}(\alpha_1,...,\alpha_k)}}
\newcommand{\pb}{\ensuremath{{\mathcal{P}}(0,\beta)}}
\newcommand{\pma}
{\ensuremath{\left\langle 
\left(\frac{1}{1-\bar{\alpha_i}{\alpha_j}}\right) \right\rangle 
}}
\newcommand{\pmc}
{\ensuremath{\left\langle 
\left(\frac{1}{1-\bar{\beta_i}{\beta_j}}\right) \right\rangle }}
\newcommand{\bda}{\ensuremath{{\mathcal{D}}(A(\mathbb{D}^2);z_1,...,z_k)}}
\newcommand{\dalp}{\ensuremath{{\mathcal{D}}(A(\mathbb{D});\alpha_1,...,\alpha_k)}}
\newcommand{\dbet}{\ensuremath{{\mathcal{D}}(A(\mathbb{D});\beta_1,...,\beta_k)}}
\newcommand{\ib}{\ensuremath{{\mathcal{D}}(A;x_1,...,x_k)}}
\newcommand{\cA}{\ensuremath{{\mathcal{A}}}}
\newcommand{\cB}{\ensuremath{{\mathcal{B}}}}
\newcommand{\cO}{\ensuremath{{\mathcal{O}}}}
\newcommand{\cH}{\ensuremath{{\mathcal{H}}}}
\newcommand{\cP}{\ensuremath{{\mathcal{P}}}}
\newcommand{\cI}{\ensuremath{{\mathcal{I}}}}
\newcommand{\cV}{\ensuremath{{\mathcal{V}}}}
\newcommand{\cM}{\ensuremath{{\mathcal{M}}}}
\newcommand{\Bd}{\ensuremath{{\mathcal{B}_{\delta}}}}
\newcommand{\D}{\ensuremath{{\mathcal{D}}}}
\newcommand{\I}{\ensuremath{{\mathcal{I}}}}
\newcommand{\Set}{\ensuremath{{\mathcal{S}}}}
\newcommand{\disk}{\ensuremath{{\mathbb{D}}}}
\newcommand{\ck}{\ensuremath{{\mathbb{C}}^k}}
\newcommand{\p}{\ensuremath{\perp}}
\newcommand{\PB}{\ensuremath{{\mathcal{P}}}}
\newcommand{\CC}{\ensuremath{\mathbb{C}}}
\newcommand{\MFB}{\ensuremath{\mathfrak{B}}}
\newcommand{\MFV}{\ensuremath{\mathfrak{V}}}
\newcommand{\MFO}{\ensuremath{\mathfrak{O}}}
\newcommand{\MFH}{\ensuremath{\mathfrak{H}}}
\newcommand{\MFP}{\ensuremath{\mathfrak{P}}}

\begin{abstract}
We compare and classify various types of Banach algebra norms on $\CC^k$
through geometric properties of their unit balls. This study is motivated by various open problems in interpolation theory and in the isometric characterization of operator algebra norms.

\end{abstract}

\section{Introduction}

Given two points $v=(v_1,...,v_k)$ and $w=(w_1,...,w_k)$ in 
$\CC^k$
we define $v\cdot w:=(v_1w_1,...,v_kw_k)$ so that $\CC^k$ is a 
commutative algebra with unit
$e=(1,...,1)$.  If in addition we have a norm $\| \cdot \|$ on
$\CC^k$ which satisfies $\|v\cdot w\| \leq \|v\| \cdot \|w\|$ for
all $v,w$ in $\CC^k$, then the pair $(\CC^k, \| \cdot \|)$
is a Banach algebra. We shall always require that $\|e\| =1.$ 
Given a norm $\| \cdot \|$ on $\CC^k$ we will refer
to the set $B:= \{ v \in \CC^k: \| v \| \leq 1 \}$ as the 
$unit$ $ball$
in $\CC^k$ with respect to the given norm. 
We will denote the set of all balls determined
by Banach algebra norms on $\CC^k$ as $\cB_k$.

In this paper we begin to compare and classify different 
properties that Banach algebra norms on $\CC^k$ can possess in 
order to provide some insight into a number of problems. One of 
our motivations is to gain a deeper understanding of 
the geometry of the balls that arise as the solutions to 
various interpolation problems. The other 
problems include an attempt to characterize the Banach algebras 
that can be represented isometrically as algebras of operators 
on a Hilbert space, questions about the computability of 
various interpolation problems for uniform algebras, and 
questions about the failure of the multi-variable von Neumann 
inequality for three or more contractions on Hilbert space.

\begin{definition}
Let $\cA$ be a unital Banach algebra. Then $\cA$
is said to {\bf satisfy von Neumann's inequality} provided that 
for all 
$a$ in $\cA$ with $\| a \| \leq 1$ we have that
$$
\| p(a) \| \leq \|p \|_{\infty}
$$
where $p$ is a polynomial and $\| p \|_{\infty}= \sup\{ |p(z)| : 
\, |z| \leq 1 \}$.
\end{definition}

Von Neumann proved that if $H$ is a Hilbert space, then the 
algebra of operators on $H, B(H)$ satisfies von 
Neumann's inequality. Consequently, any unital Banach algebra 
that can be represented isometrically as  an algebra of 
operators on a Hilbert space necessarily satisfies von 
Neumann's inequality. It is still unknown if the converse of 
this last statement holds. That is, can every unital Banach 
algebra that satisfies von Neumann's inequality be represented 
isometrically as an algebra of operators on a Hilbert space?
It has been conjectured that the answer to this question is 
affirmative.
For this reason we are interested in understanding the Banach 
algebra norms on $\CC^k$ that satisfy von Neumann's inequality.

\begin{definition}
The set of all balls determined by Banach algebra norms on 
$\CC^k$ satisfying
von Neumann's inequality will be denoted $\cV_k$.
The set of all balls determined by Banach algebra norms on 
$\CC^k$ such that the resulting Banach algebra has a unital 
isometric representation as an algebra of operators on a 
Hilbert space will be denoted $\cO_k.$
\end{definition}

By von Neumann's result, 
we have that $\cO_k \subseteq \cV_k \subseteq \cB_k$.
If the answer to the above question is affirmative, then it 
must be the case that $\cO_k = \cV_k.$ Conversely, if we were 
able to show, for some $k$, that these two sets are not equal, 
then that would provide a counterexample to the above 
conjecture.

The next subset of $\cB_k$ that we wish to examine is motivated
by general interpolation theory and by attempts to obtain
multi-variable generalizations of
von~Neumann's inequality.

\begin{definition}
A unital, commutative Banach algebra, $\cA$, is said to {\bf 
satisfy the multi-variable von Neumann inequality} provided 
that for every $n$ and for every set of $n$ elements $\{a_1, 
\ldots, a_n \}$ from the unit ball of $\cA$ and for every 
polynomial $p$ in $n$ variables,
we have that
$$\|p(a_1, \ldots, a_n)\| \le \|p\|,$$
where $\|p\| = sup\{|p(z_1, \ldots, z_n)| : |z_i|\le 1, 1 \le i 
\le n \}.$
A subset $B$ of $\CC^k$ is called {\bf hyperconvex} if it is 
the unit ball of a Banach algebra norm on $\CC^k$ that 
satisfies the multi-variable von Neumann inequality and we let 
$\cH_k$ denote the collection of all such balls.
\end{definition}

Let $X$ be a compact Hausdorff space and let $C(X)$ denote the
continuous complex-valued functions on $X$. We call $A \subseteq
C(X)$ a $uniform$ $algebra$ provided that $A$ is uniformly 
closed,
contains the identity, and separates points in $X$.
An algebra that is isometrically isomorphic to the quotient of 
a 
uniform algebra
is called in the literature, we are sorry to say, an {\em IQ 
algebra}. 
Davie\cite{Dav} proved that a unital, commutative Banach 
algebra with a unit of norm 1,
is an IQ algebra if and only if it satisfies the multi-variable 
von~Neumann inequality.

The term {\em hyperconvex} is due to Cole and Wermer\cite{CW1}, 
who introduced this concept because of the relationship between 
these sets, isometric 
quotients of uniform algebras and questions in interpolation 
theory for uniform algebras. We explain these connections below.

Ando\cite{And} proved a 2 variable analogue of von Neumann's 
inequality, namely, that every commuting pair of contractions
on a Hilbert space satisfies the 2-variable von Neumann 
inequality. Since then many mathematicians have given examples 
to show that commuting triples of contractions, even on a 
finite dimensional Hilbert 
space need not satisfy the 3-variable von Neumann inequality.
Thus, a commuting algebra of operators on a Hilbert space need 
not satisfy the multi-variable von Neumann inequality, i.e., 
need not be an IQ algebra.
However, these examples generally consisted of commuting sets
of nilpotent operators and examples of commuting sets of 
diagonalizable
operators have been much harder to come by.
In particular, it is still unknown if a simultaneously 
diagonalizable set of three
commuting contractions 
on a three dimensional Hilbert space must satisfy
the multi-variable von~Neumann inequality. 
Lotto and Steger\cite{LS} constructed three commuting, 
diagonalizable contractions on a five
dimensional Hilbert space that failed to satisfy von~Neumann's 
inequality.
Recently, Holbrook \cite{Hol} constructed three commuting, 
diagonalizable contractions on a four
dimensional Hilbert space that failed to satisfy von~Neumann's inequality.
We shall see that this problem  is closely related to 
determining whether or
not $\cO_k=\cH_k$.

Cole,  
Lewis, and Wermer
in \cite{CLW} make the following definition.

\begin{definition}
Let $A \subseteq C(X)$ be a uniform algebra and let $x_1,...,x_k 
\in X$. The {\bf interpolation
body} associated with $A$ and $x_1,...,x_k$, denoted as $\ib$, 
is 
defined as 
$$
\{ (w_1,...,w_k): \forall \, \epsilon >0 \, \exists \, f \in A 
\, \ni \,
\|f \|_{\infty} \leq 1+\epsilon \, \mathrm{ and } \, f(x_i)=w_i 
\}
$$
\end{definition}

The connection between interpolation bodies and hyperconvex 
sets is due to Cole and Wermer \cite{CW1}, it uses the 
characterization of IQ algebras due 
to S. Davie \cite{Dav}.

Let $I_x$ denote the ideal of functions vanishing at
the points $x_1,...,x_k$ in $X$. 
Since $A$ separates points, there exist functions
$f_1,...,f_k$ in $A$ with $f_i(x_j)=\delta_{ij}$ . Thus,
in $A/I_x$ we have that $F_i=[f_i + I_x]$ is a set of 
$k$-commuting
idempotents,
satisfying $F_iF_j=\delta_{ij}F_j$, $F_1 + \hdots + F_k=I$, and 
which span
$A/I_x$. 

The following is immediate.
\begin{proposition} \label{interpbody}
Let $A \subseteq C(X)$ be a uniform algebra and let 
$x_1,...,x_k 
\in X$.
Then
$$
\D(A;x_1,...,x_k)= 
\{ (w_1,...,w_k) \in \CC^k: \|w_1F_1 + \cdots +w_kF_k \| \leq 1 
\}.
$$
\end{proposition}

Hence, a point $(w_1,...,w_k)$ in $\ck$
belongs to $\ib$ if and only if
$\| w_1F_1 + \hdots + w_kF_k \| \leq 1$.
Thus, an interpolation body $\ib$ is a natural 
coordinatization of the
closed unit ball of the quotient algebra
$A/I_x$.
The isomorphism between $A/I_x$ and $\ck$
defined by sending a coset $[f+I_x]$ to $(f(x_1), \ldots , 
f(x_k))$
endows $\ck$ with a Banach algebra norm for which $\ib$ is the 
closed unit
ball.

Combining Proposition~1.5 with Davie's characterization of IQ algebras yields the following result of Cole-Wermer.

\begin{proposition}(\cite{CW1}, Theorem~5.2)
Let $A \subseteq C(X)$ be a uniform algebra and let 
$x_1,...,x_k 
\in X$, then $\ib$ is a hyperconvex subset of $\CC^k$.  
Conversely, if $B$ is any hyperconvex subset of $\CC^k$, then 
there is a uniform algebra $A \subseteq C(X)$ and points,
$x_1, \ldots, x_k \in X$, such that $B = \ib.$
\end{proposition}

Thus, the hyperconvexity condition gives an abstract 
characterization of interpolation bodies. Unfortunately, given a ball $B$
in $\CC^k$ it is very difficult to determine if it is hyperconvex. One
of the main difficulties is that to apply the hyperconvexity
condition one needs analogues of the Nevanlinna-Pick interpolation theorem on
polydisks, which is still an open problem for more than two variables.

By a result of Cole \cite{BD}, page 272,
$A/I_x$ is an operator 
algebra. 
Thus, we have that $\cH_k \subseteq \cO_k$ for all $k$.

Another collection of Banach algebra balls  that we wish to 
consider
is a 
well understood subset of $\cH_k$. 
When the uniform algebra $A$ is the disk algebra, denoted
$A(\mathbb{D})$, the set $\dalp$ is referred to as a {\em Pick 
body} 
and will be denoted as 
$\pa$. By Pick's theorem \cite{Pic} , if $\{\alpha_1, \ldots, 
\alpha_k\}$
is a subset of the open unit disk, then a point $(w_1,...,w_k)$ 
belongs
to $\pa$ if and only if the matrix $\left( \frac{1-\bar{w_i}w_j}{
1-\bar{\alpha_i}\alpha_j} \right)$ is positive semi-definite. 
When these
points are all on the unit circle, then $\pa$ can be shown to be 
the closed
$k$-polydisk.
The set of balls in $\CC^k$ determined by interpolation bodies 
of the form $\pa$
will be denoted as $\cP_k$. 

Thus, to summarize, we know that for any $k$ in 
$\mathbb{N}$ we have 
that
$$
\cP_k \subseteq \cH_k \subseteq \cO_k \subseteq \cV_k 
\subseteq \cB_k.
$$
Moreover, we have argued that proving equality between various 
pairs 
of these subsets will yield
positive results on some open problems, while showing that 
certain of these subsets are not
equal will yield counterexamples to certain open questions.

Thus, we are lead to study various properties of the balls belonging to these five families, to attempt to obtain more manageable characterizations of these five families, and to find means of generating balls that belong to these various families, in an attempt to distinguish between these five families.
We will see that in some senses, focusing on the properties of the balls is an efficient way to produce examples of interpolation sets and operator algebras without actually needing to construct the uniform algebras or operators.
 
Here is a summary of what we can show about the relationships between
these sets:
\begin{eqnarray*}
k=2 & \hspace{.5cm} & \cP_2 = \cH_2 = \cO_2 = \cV_2 
\subsetneq \cB_2\\
k=3 & \hspace{.5cm} & \cP_k \subsetneq \cH_k \subseteq^?
\cO_k \subseteq^? \cV_k \subsetneq \cB_k\\
k \geq 4 & \hspace{.5cm} & \cP_k \subsetneq \cH_k \subsetneq
\cO_k \subseteq^? \cV_k \subsetneq \cB_k\\
\end{eqnarray*}
Thus, for $k=2$, we will show that $\cP_2=\cH_2=\cO_2=\cV_2$ and 
provide an example
of a Banach algebra norm in 2 dimensions that fails to satisfy 
von~Neumann's inequality.
For $k \geq 4$, we are able to show that the sets $\cP_k,
\cH_k,
\cO_k,$ and $\cB_k$ are distinct. However, for $k=3$, the 
picture is not complete and we are still unable to resolve 
whether or not $\cO_k =\cV_k$ for all
$k \geq 3$.

Some of these results are restatements of earlier results.
In \cite{CLW} (Theorem 5), Cole, Lewis, and Wermer show
that every hyperconvex set $Y$ in $\mathbb{C}^2$ is a Pick
body, i.e., that $\cH_2 = \cP_2$. 
In fact, we shall see that their proof actually implies
$\cP_2=\cV_2$.   They also show that for $k>2$,
$\cH_k \neq \cP_k$. Similarly, the fact that $\cH_k \ne \cO_k$ for $k \ge 4,$ is a restatement of a result of Holbrook\cite{Hol}.

There is one final concept that we shall study.
Recall that a set is called {\em semi-algebraic} if it is 
defined by a finite set of
polynomial inequalities.
In a certain sense these are the most definable or computable 
subsets of $\ck,$
since determining whether or not a point belongs to 
such a set involves 
only finitely many 
algebraic operations, provided that the polynomials can be found.

It is easily seen that all the Pick bodies are semi-algebraic 
sets.
In \cite{CW2}, Cole and Wermer show that 
$\D(A(\mathbb{D}^2);z_1,...,z_k)$
is a semi-algebraic set. 
If $A$ denotes the uniform algebra of all bounded, analytic 
functions on an annulus, then results
of \cite{FV} show that $\ib$ is semi-algebraic for any choice of 
points.

It is natural to wonder if every interpolation body is 
semi-algebraic.
If interpolation bodies are known a priori to be semi-algebraic 
sets,
then we would know that solving the corresponding interpolation 
problem can  always be reduced
to a finite set of algebraic operations.
We conjecture that for a sufficiently pathological uniform 
algebra, there will exist
interpolation bodies that are not semi-algebraic.
The most interesting question is whether or not an interpolation 
body exists for some natural
uniform algebra that is not semi-algebraic.
However, since every two point interpolation body is a Pick 
body, we see that these
interpolation problems are in some sense computable.

Although 
we 
have been unable to produce an example of an interpolation body 
that is not
semi-algebraic, we will provide an example of 
a ball in $\cO_3$ which we conjecture is not a semi-algebraic set.
Thus, if the before mentioned conjecture is true, then either $\cH_3 \ne \cO_3$ or there exist 3 point 
interpolation problems for which the
interpolation body is not semi-algebraic.

\section{Characterization and Generation of Sets in These Families}

In this section we present some basic results about the various 
classes of norms on $\ck$ that
we have introduced. We begin by gathering together some facts 
about $\cB_k.$

Recall that a set $B \subseteq \ck$ is the closed unit ball of 
some norm on $\ck$ if and only
if $B$ is closed, bounded, absorbing and absolutely convex. We 
shall refer to such a set as
a {\bf ball.}
Given a finite collection of balls, it is easily seen that their 
intersection will 
be closed, bounded, absorbing and absolutely
convex. Thus, the intersection of a  finite collection of balls 
is again a ball. However, for
arbitrary collections of balls, their intersection will be 
closed, bounded and absolutely convex, but
not necessarily absorbing.
Thus, the intersection will be a ball if and only if
it is absorbing. We shall show a similar result holds for the 
various classes of norms that
we wish to consider.

Given a set $B \subseteq \ck$ we set $B \cdot B = \{ v\cdot w : 
v,w \in B \}.$

\begin{proposition}
Let $B \subseteq \ck.$ Then $B \in \cB_k$ if and only if $B$ is 
a ball, $e \in B$ and $B\cdot B
 \subseteq B.$ Consequently, the intersection of a  finite 
collection of sets in $\cB_k$ is
 again in $\cB_k$, but the intersection of an arbitrary 
collection of sets
 in $\cB_k$ is in $\cB_k$ if and only if the intersection is 
absorbing.
 \end{proposition}
 
Note that the fact that $B \subseteq \ck$ is bounded and $B\cdot B \subseteq 
B$ implies that $B$ must be a
subset of the closed unit polydisk. 

\begin{definition}
Given any set $D \subseteq \ck$ contained in the closed unit polydisk, we 
let $\cB(D)$
denote the intersection of all elements of $\cB_k$ containing
$D$.
\end{definition}

\begin{definition} Let $D \subseteq \ck$. We say that $D$ is 
{\bf separating} provided that for every $i \ne j$ there exists 
$(w_1, \ldots, w_k) \in D$ such that $w_i \ne w_j.$
\end{definition}

Note that $D= \{v_1, \ldots, v_n \} \subseteq \ck$ is 
separating where $v_i = (w_{i,1}, \ldots, w_{i,k})$, if and 
only if setting $x_j = (w_{1,j}, \ldots, w_{n,j})$ defines $k$ 
distinct points in $\mathbb{C}^n.$
 
\begin{proposition} \label{intersection}
Let $D \subseteq \ck$ be a separating set contained in the 
closed 
unit polydisk. Then 
$\cB(D) \in \cB_k$ and if $D \subseteq B_1$ with $B_1 \in \cB_k$ 
then $\cB(D) \subseteq B_1.$
\end{proposition}
\begin{proof}
Clearly, if $B_1 \in \cB_k$ contains $D$, then $\cB(D) 
\subseteq B_1.$ So it remains to show that $\cB(D) \in \cB_k,$ 
which 
by the above is equivalent to proving that $\cB(D)$ is 
absorbing. 
Since $D$ is separating, $\cB(D)$ is separating. Fix $i \ne 
j$, then by taking an absolute convex combination with $e$, we 
have an element of $\cB(D)$ that is a strictly positive constant 
in 
the i-th coordinate and 0 in the j-th coordinate. Fixing $i$, 
choosing one such vector for every $j \ne i$ and taking the 
product of these vectors yields a vector in $\cB(D)$ that is 
strictly positive in the i-th coordinate and is 0 in the 
remaining coordinates. 

Repeating for each $i$, we see that $\cB(D)$ contains a positive 
multiple
of each basis vector. Since the absolutely convex hull 
of the set of  such vectors will contain an open neighborhood of 
0, we have that 
$\cB(D)$ is absorbing.
\end{proof}

When $D$ is separating and a subset of the closed unit polydisk, 
we shall refer to $\cB(D)$ 
as the {\em Banach algebra ball generated by D.}
Later, we shall prove that $\cB(D)$ is a ball if and only if $D$ 
is separating.

\begin{example}
Let $D = \{ e_1,e_2 \} \subseteq \CC^2$ where $e_1=(1,0), 
e_2=(0,1)$ which is a separating set.
Clearly, $\cB(D)$ must contain the absolutely convex hull of 
$e_1, e_2$ and $e$. But 
we claim that this latter set
is closed under products and hence is in $\cB_2.$ To see this 
last claim, note that if
$w=w_1e_1+w_2e_2+w_3e_3$ and $z=z_1e_1+z_2e_2+z_3e_3$ with 
$|w_1|+|w_2|+|w_3| \le 1$ and
$|z_1|+|z_2|+|z_3| \le 1$, then 
$$w \cdot z = (w_1z_1+ w_1z_3 +w_3z_1)e_1 + (w_2z_2+w_2z_3+ 
w_3z_2)e_2 + (w_3z_3)e$$
and the sum of the absolute values of the coefficients of $w 
\cdot z$ are easily seen to be less than 1.
\end{example} 

We now turn our attention to $\cV_k.$ Given a set $V \subset 
\ck$ and a function $f$, we set
$f(V) = \{ (f(v_1), \ldots, f(v_k)) : (v_1, \ldots, v_k) \in V 
\}.$

\begin{proposition}
Let $V \subseteq \ck.$ Then $V \in \cV_k$ if and only if $V \in 
\cB_k$ and $f(V) \subseteq V$
for every $f \in A(\disk)$ with $\|f\| \le 1.$
Consequently, the intersection of a  finite collection of sets 
in $\cV_k$ is again in $\cV_k$,
but the intersection of an arbitrary collection of sets in 
$\cV_k$ is in $\cV_k$
if and only if the intersection is absorbing.
\end{proposition}

\begin{definition} Given any set $D \subseteq \ck$ contained in 
the closed unit polydisk,
we let $\cV(D)$ denote the intersection of all elements of 
$\cV_k$ that contain $D.$
\end{definition}

\begin{proposition}
Let $D \subseteq \ck$ be a separating set contained in the 
closed unit 
polydisk. Then 
$\cV(D) \in \cV_k$ and if $D \subseteq V_1$ with $V_1 \in \cV_k$ 
 then $\cV(D) \subseteq V_1.$
\end{proposition}

\begin{proof}
As before it will be enough to prove that $\cV(D)$ is absorbing. 
But note that $\cB(D) \subseteq \cV(D)$ and that $\cB(D)$ is 
absorbing. Hence, $\cV(D)$ is 
absorbing and so 
the result follows.
\end{proof}

When $D$ is separating and contained in the closed unit polydisk,
we shall refer to $\cV(D)$ as the {\em von~Neumann algebra ball 
generated by D.}
We shall also see that $\cV(D)$ is a ball if and only if $D$ is 
separating.

We let $\varphi_a(z) = \frac{z-a}{1- \bar{a}z}$ denote the 
elementary Mobius map. Using these maps
it is possible to reduce problems about $\cV_k$ by one-dimension.

Given a set $V \subset \ck,$ we set $\widehat{V} = \{(v_2, 
\ldots, v_k) \in \CC^{k-1}: (0, v_2, \ldots, v_k) \in V\}.$

\begin{proposition}
Let $V_1,V_2 \in \cV_k$. Then $V_1=V_2$ if and only if 
$\widehat{V_1} = \widehat{V_2}.$
\end{proposition}

\begin{proof}
Assume that $\widehat{V_1} =\widehat{V_2}$.
If $(a, v_2, \ldots, v_k) \in V_1$, then $(0, \varphi_a(v_2), 
\ldots, \varphi_a(v_k)) \in V_1$ and
hence in $V_2$. Applying $\varphi_{-a}$ we find that $(a, v_2, 
\ldots, v_k) \in V_2$.
Thus, $V_1 \subseteq V_2,$ but by the symmetry of the argument, 
$V_1=V_2.$

The other implication is obvious.
\end{proof}

\begin{theorem}
We have that $\cP_2 = \cH_2 = \cO_2 = \cV_2 \ne \cB_2.$
\end{theorem}

\begin{proof}
Let $V \in V_2.$ Since $V$ is a ball $\widehat{V} = \{ z \in 
\CC : |z| \le r \}$ for
some $r, 0< r \le 1.$ Hence, $\widehat{V} = 
\widehat{\cP(0,r)}$ and so $V= \cP(0,r).$ Thus, $\cV_2 
\subseteq \cP_2,$
and the equalities  follow.

To see that $\cV_2 \ne \cB_2$, we consider the Banach algebra 
ball
$\cB(D)$ for $D = \{ e_1, e_2 \}$ constructed in the above 
example.
Assume $\cB(D) \in \cV_2$, then we see that $\widehat{\cB(D)}$ 
is the closed unit
disk and so necessarily, $\cP(0,1) = \cB(D).$ However, it is 
easily shown that $\cP(0,1)$
is the closed unit bidisk and so $\cP(0,1) \ne \cB(D).$ Thus, 
$\cB(D)$ is not in $\cV_2.$ 
\end{proof}

The proof that $\cP_2 = \cV_2$ is essentially the one given by 
Cole and Wermer 
to prove that $\cP_2=\cH_2.$

\begin{problem}
Find necessary and sufficient conditions on a set $B \subseteq 
\CC^{k-1}$ such that
$B=\widehat{V}$ for some $V \in \cV_k.$
\end{problem}

Given $B \subseteq \CC^{k-1},$ we set
$$\widetilde{B}= \{(a, v_2, \ldots, v_k) : (\varphi_a(v_2), \ldots, \varphi_a(v_k)) \in B\},
$$
so the problem is to find conditions on $B$ that guarantees $\widetilde{B} 
\in \cV_k.$

The following gives some conditions that $B$ must necessarily 
satisfy.

\begin{proposition}
If $B=\widehat{V}$ for some $V \in \cV_k$, then
\begin{enumerate}
\item $B$ is a closed subset of the closed polydisk,
\item $B$ is a ball in $\CC^{k-1}$,
\item $B \cdot B \subseteq B,$
\item for any $f \in A(\disk)$ with $\|f\| \le 1$ and $f(0)=0$, 
we have $f(B) \subseteq B.$
\end{enumerate}
\end{proposition}

Note that these conditions imply that $B$ is the unit ball of a 
Banach algebra norm on $\CC^k$, but one for
which, possibly, $\|e\| > 1.$

Using the Nevanlinna-Pick theorem it is possible to replace the 
last condition on $B$ by a matrix positivity condition that 
makes no reference to analytic functions. In fact, condition 4) 
is equivalent to requiring that whenever, $(v_1, \ldots, 
v_{k-1}) \in B$ and
$$\begin{pmatrix} \frac{\bar{v_i}v_j - \bar{w_i}w_j}{1- 
\bar{v_i}v_j} \end{pmatrix}$$ is 
positive semi-definite, then $(w_1, \dots, w_{k-1}) \in B.$

The following example shows that these conditions on $B$ are not 
sufficient and illustrates some of the difficulties
in determining even the elements of $\cV_3.$

\begin{example}
Fix $0 < r \le 1$ and let $B= \{(w_1, w_2) : |w_1| + |w_2| \le r 
\}.$
We will show that $B$ satisifes the conditions of the above 
proposition, but that $ \widetilde{B}=V$
is not convex.
\end{example}
Clearly, $B$ satisifes the first two conditions. If 
$w=(w_1,w_2)$ and $z=(z_1,z_2)$ are in $B$,
then $|w_1z_1| + |w_2z_2| \le |w_1| + |w_2| \le r$ and so 
$w\cdot z \in B.$
Also, if $f \in A(\disk)$, with $\|f\| \le 1$ and $f(0)=0$, 
then $f(z)=zg(z)$ with $\|g\| \le 1.$
>From this it is easily seen that $|g(w_1)|+ |g(w_2)| \le r.$

Note that for any, $|a| \le 1, (a,b,a) \in V$ if and only if 
$|\varphi_a(b)| \le r.$
Similarly, $(-a,-a,c) \in V$ if and only if $|\varphi_{-a}(c)| 
\le r.$
If $V$ was convex, then $1/2[ (a,b,a) + (-a,-a,c)] \in V$ and 
consequently,
$|b-a| + |c+a| \le 2r.$ However, for $0<a<r,$ if we choose $b,c$ satisfying $\varphi_a(b)= -r, \varphi_a(c)= +r,$ then this latter inequality fails.

Similar calculations show that for $B= \{(w_1,w_2) : |w_1|^2 + 
|w_2|^2 \le r^2 \}$ one also has that
$ \widetilde{B}$ is not convex.

Note that in both of these examples the set $B$ is {\em 
circled}, that is $(w_1,w_2) \in B$
implies that $(e^{i\theta_1}w_1, e^{i\theta_2}w_2) \in B.$

\begin{problem} Find necessary and sufficient conditions on a 
circled set $B \subseteq \CC^2$ such that
$B=\widehat{V}$ for some $V \in \cV_3.$
\end{problem}

We now examine $\cO_k.$ There is a dual object for sets in 
$\cO_k$, called a Schur ideal,
which we will find makes it much easier to determine if a set 
belongs to $\cO_k$ than to
$\cV_k$ or $\cH_k.$
The following gives a convenient description of elements of 
$\cO_k.$

\begin{definition}
Let $E_1,...,E_k$ be bounded operators on a Hilbert space 
$\mathbf{H}$
satisfying
\begin{enumerate} [(a)]
\item $E_iE_j = \delta_{ij}E_i$ and
\item $E_1 + \cdots + E_k =I$.
\end{enumerate}
Then the algebra $\cA :=\mathrm{span} \{E_i\}_{i=1}^k$ is 
called a {\bf $k$-idempotent operator algebra}.
\end{definition}

Each $k$-idempotent operator algebra determines a ball in 
$\CC^k$ in
the following way.

\begin{proposition} \label{operatorballs}
Let $\cA$ be a $k$-idempotent operator algebra. Then the set
$$
\D(\cA ):= \{(w_1,...,w_k)\in \CC^k: \| w_1E_1 +\cdots + w_kE_k 
\| \leq 1 \}
$$
has the following properties:
\begin{enumerate} [(a)]
\item $\D(\cA ) \in \cO_k$,
\item $\| v \|_{\D(\cA )}:=\mathrm{inf} \{ t : v \in t\D(\cA ) 
\}$ 
is a 
Banach algebra norm on $\CC^k$, and 
\item $\| \cdot \|_{\D(\cA )}$ satisfies von Neumann's 
inequality.
\end{enumerate}
Moreover, $B \in \cO_k$ if and only if $B= \D(\cA)$ for some 
k-idempotent operator algebra $\cA.$
\end{proposition}

\begin{proposition}
The intersection of a finite collection of sets in $\cO_k$ is 
again in $\cO_k$, but the intersection
of an arbitrary collection of sets in $\cO_k$ is in $\cO_k$ if 
and only
if the intersection is absorbing.
\end{proposition}
\begin{proof}
Let $B_{\alpha}$ be a collection of sets in $\cO_k$ and for each 
$\alpha$ choose a
$k$-idempotent operator algebra $\cA_{\alpha}$  acting on a 
HIlbert space $\cH_{\alpha}$ and
generated by $\{E_i^{\alpha} \}$ such that
$B_{\alpha} = \D(\cA_{\alpha})$.
 
If the operators, $E_i = \oplus_{\alpha} E_i^{\alpha}$ acting 
on the direct sum of the Hilbert
spaces are bounded, then they will clearly generate a 
$k$-idempotent operator algebra $\cA$ with
$\D(\cA)$ equal to the intersection of the collection. 
Thus, all that remains is to note that the boundedness of the 
set of operators $E_i$ is equivalent to
the fact that the intersection is absorbing.
\end{proof}

\begin{definition} Given $D \subseteq \ck$ contained in the 
closed unit polydisk. we let
$\cO(D)$ denote the intersection of all set in $\cO_k$ that 
contain $D.$
\end{definition}

\begin{proposition} Let $D \subseteq \ck$ be a subset of the 
closed polydisk. If $D$ is separating, then $\cO(D) \in \cO_k.$ 
\end{proposition}
\begin{proof}
Again, it is enough to prove that $\cO(D)$ is absorbing, but 
since
$\cB(D) \subseteq \cO(D)$, the result follows.
\end{proof}

When $D$ is separating and a subset of the closed unit polydisk, 
then we refer to
$\cO(D)$ as the {\em operator algebra ball generated by D.}

In \cite{Pau1}, the first author developed a dual object for 
operator norms,
called
$Schur$ $ideals$. These will give us an alternate description of 
$\cO(D)$.

Let $\I$ be
a subset of $M_k^+$, the positive semi-definite $k \times k$ 
matrices satisfying:
\begin{enumerate} [i)]
\item $P,Q \in \I$ then $P+Q \in \I$
\item $P \in \I, Q$ positive, then $P \ast Q \in \I$
\end{enumerate}
where $(p_{ij})\ast (q_{ij})=( p_{ij}q_{ij})$, the Schur 
product, then $\I$ is called a 
{\bf Schur ideal}. Note
that the set of positive matrices is closed under $+$ and 
$\ast$, and
that the matrix of all $1's$ acts as a unit for $\ast$. Given a 
set 
$\Set$ of positive matrices we let $\langle \Set \rangle$ denote 
the
Schur ideal that it generates. 

Given a set $\Set$ of
positive $k\times k$ matrices and a subset $\D$ of 
the closed $k$-polydisk containing $0$, the first author in 
\cite{Pau1} 
defines
$$
\Set^{\p}:= \{(w_1,...w_k)\in \CC^k:((1-\bar{w_i}w_j)p_{ij}) 
\geq 0
\text{ for all } (p_{ij}) \in \Set \} 
$$
and
$$
\D^{\p}:= \{ (p_{ij}) : ((1-\bar{w_i}w_j)p_{ij}) \geq 0 \text{ 
for all }
(w_1,...,w_k) \in \D \}.
$$

Note that $\Set^{\p} = \langle \Set \rangle^{\p}.$

In \cite{Pau1} it is shown that,
if $\D \in \cO_k$, then 
$\D^{\p}$ is a Schur ideal with the following two properties:
\begin{enumerate} [a)]
\item for each $i=1,...k$ there exists $P \in \D^{\p}$ such that 
$p_{ii} \not= 0$ and
\item there exists a $\delta >0$ such that for all $P \in 
\D^{\p}$, we have
that $P \geq \delta^2 \mathrm{Diag}(P).$
\end{enumerate}
A Schur ideal satisfying $a)$ is said to be {\bf non-trivial}. A 
Schur
ideal satisfying $b)$ is said to be {\bf bounded}.

Given an arbitrary Schur ideal (which is $non-trivial$ and 
$bounded$) contained
in $M_k^+$, say $\I$, one can construct an ``affiliated 
''$k$-idempotent 
operator algebra, say $\cA_{\I}$. The following proposition, 
which is a non-matricial version
of Theorem 3.2 in \cite{Pau2}
provides a way to
compute the ball $\D(\cA_{\I})$ in terms of the ``perp'' of the 
Schur ideal $\I$.

\begin{proposition} \label{schuralgebra}
Let $\I$ be a non-trivial and bounded Schur ideal contained in 
the
positive $k\times k$ matrices. Then there exists
a $k$-idempotent operator algebra $\cA_{\I}$ such that 
$$
\D(\cA_{\I}) = \I^{\p}.
$$
\end{proposition}

\begin{proof} First we will construct the $k$-idempotent 
operator 
algebra $\cA_{\I}$. Let $\I^{-1}$ denote the set of invertible
elements in $\I$ and $E_{ii}$ denote the canonical matrix units. 
Note that
since  $\I$ is non-trivial the set of invertible elements 
$\I^{-1}$ are
dense in $\I$. This can be seen by observing the identity matrix 
$I$ belongs
to $\I$ and we have that for any $Q \in \I$, $(Q+\epsilon I) \in 
\I$ and is
invertible.
Now set $E_i := \bigoplus_{Q \in \I^{-1}}
Q^{1/2}E_{ii}Q^{-1/2}$ for $i=1,...,k$. These are operators 
acting on 
$\bigoplus_{P \in \I^{-1}} M_k$ and we let
$$
\cA_{\I} : = \mathrm{span} \{ E_1,...,E_k\}.
$$
Now since $\I$ is bounded, by definition we have that there 
exists
a $\delta >0$ such that $(q_{ij}) \geq 
\delta^2\mathrm{Diag}(q_{ii})$ for all 
$(q_{ij}) \in \I^{-1}$. Thus, $\delta^2q_{mm}E_{mm} \leq 
(q_{ij})$ for
each $1 \leq m \leq k$ and for each $(q_{ij}) \in \I^{-1}$. 
We have that $E_i^*E_i =
\bigoplus_{Q \in \I^{-1}} Q^{-1/2}q_{ii}E_{ii} Q^{-1/2} \leq
\delta^{-2} \bigoplus_{Q \in \I^{-1}} I_{k \times k}$ and
we have that $\|E_i\| \leq \delta^{-1}$. Thus each $E_i$ is 
bounded and
$\cA_{\I}$ is a $k$-idempotent operator algebra.

We now have the following logical equivalences:
\begin{eqnarray*}
(w_1,...,w_k) \in \D(\cA_{\I}) 
& \Longleftrightarrow & \left\| \sum\limits_{i=1}^k w_iE_i 
\right\| \leq 1 \\
& \Longleftrightarrow &
\left\| \bigoplus\limits_{Q \in \I^{-1}} 
Q^{1/2}\mathrm{Diag}(w_i) Q^{-1/2} \right\| \leq 1 \\
& \Longleftrightarrow & 
((1- \bar{w_i}w_j)q_{ij}) \leq 0 \text{ for all } Q \in \I^{-1}\\
& \Longleftrightarrow & (w_1,...,w_k) \in \I^{\p}
\end{eqnarray*}
\end{proof}

The following is often a useful way to generate non-trivial, 
bounded Schur ideals.

\begin{proposition} Let $\Set$ be a uniformly bounded collection 
of positive, invertible matrices whose inverses are
 are also uniformly bounded. Then $\langle \Set \rangle$ is a 
non-trivial, bounded Schur ideal and consequently,
$\Set^{\p} \in \cO_k$.
\end{proposition}

\begin{proof}
Since each matrix is invertible, non-triviality follows. The 
fact that the matrices and their inverses are uniformly bounded
implies that there are positive constants $c$ and $d$ such that 
$cI \le P \le dI$ for all $P \in \Set.$
Hence, $Diag(P) \le dI \le c^{-1}dP$ for any $P \in \Set.$
Thus, there is a $\delta > 0 $ such that $P \ge \delta Diag(P)$ 
for all $P \in \Set.$
Now if $P_1, \ldots, P_m \in \Set$ and $Q_1, \ldots Q_m \in 
M_k^+$, then
\begin{eqnarray}
\nonumber
P_1*Q_1 + \cdots +P_m*Q_m & \ge & \delta Diag(P_1)*Q_1+ \cdots +\delta Diag(P_m)*Q_m \\
\nonumber
& = & \delta Diag(P_1*Q_1+ \cdots +P_m*Q_m).
\end{eqnarray}
Since $P_1*Q_1+ \cdots +P_m*Q_m$ is a typical element of the 
Schur ideal generated by $\Set$, the result follows.
\end{proof}
 
These results allow us to give a complete characterization of 
elements of $\cO_k.$

\begin{theorem}
Let $D \subseteq \ck$ be a subset of the closed polydisk. Then 
$D \in \cO_k$ if and only if $D = D^{\p \p}.$
Consequently, if $O \in \cO_k$ and $D \subseteq O$, then $D^{\p 
\p} \subseteq O$. 
\end{theorem}

\begin{proof}
The first statement is Theorem 5.8 in \cite{Pau1}.
If $D \subseteq O$, then $O^{\p} \subseteq D^{\p}$ and hence 
$D^{\p \p} \subseteq O^{\p \p} =O.$
\end{proof}

\begin{proposition} Let $D \subset \ck$ be a subset of the 
closed polydisk.
Then $D^{\p \p} = \cO(D).$
\end{proposition}
\begin{proof}
By the above we have that $D^{\p \p}$ is contained in the 
intersection.

Assume that $w= (w_1, \ldots, w_k)$ in the polydisk is not in 
$D^{\p \p}$. 
Then there exists $P=(p_{i,j}) \in D^{\p}$
such that $((1-\bar{w_i}w_j)p_{ij})$ is not positive 
semi-definite.
By continuity, we may replace P by P plus a small positive 
multiple of the identity matrix and the above matrix will
still not be positive definite.
Since the identity matrix is in $D^{\p}$, this new matrix wil be 
in $D^{\p}$.
Thus, we may assume that P is invertible. Now as in the proof of 
Proposition \ref{schuralgebra}, if we let
$\cA$ be the $k$-idempotent operator algebra generated by $\{ 
P^{1/2}E_{ii}P^{-1/2} \}$,
then $\D(\cA) \in \cO_k, D \subseteq \D(\cA)$ and $w$ will not 
be in $\D(\cA).$
\end{proof}

By the above results to determine if $\cV_k =\cO_k$, it is 
enough to decide whether or not $V^{\p \p} =V$
for every $V \in \cV_k.$

Although the above results make it relatively easy to produce 
sets in $\cO_k$ we have no clear way to determine
whether or not they are hyperconvex.

\begin{problem}
Find necessary and sufficient conditions on a Schur ideal $\cI$ 
so that $\cI^{\p}$ is hyperconvex.
\end{problem}

The following example illustrates the difficulty.

\begin{problem} Let $a,c > 0$ and set
$$P_{a,c}=
\left(
\begin{array}{ccc}
1 & 1 & 1 \\
1 & a+1 & 1 \\
1 & 1 & c+1 \\
\end{array}
\right).
$$
Then the set $\{P_{a,c} \}^{\p} \in \cO_3$, but is it 
hyperconvex?
\end{problem}

We would like
to note that one other reason that we are interested in studying 
$k$-idempotent operator algebra balls
is that often one wants to study interpolation for some operator 
algebra $A$ of 
functions on a set $X$, that is not a uniform algebra. For 
example, the algebra $\cM(X)$
of functions that act as multipliers on some reproducing kernel 
Hilbert space $\cH$  on $X$  
is a subalgebra of $B(\cH)$. This algebra equipped with the 
operator norm is sometimes, but not always,
a uniform algebra. One still wishes to study the corresponding 
interpolation bodies,
$\{(f(x_1), \dots, f(x_k)) : \|M_f\| \le 1 \}$ in this
situation. Since it is known by a result of \cite{BRS} that 
$A/I_x$ is an operator algebra 
generated by
$k$-commuting idempotents as above,
then these more general interpolation bodies will be in $\cO_k$ 
but not necessarily in $\cH_k$. 

We now turn our attention to hyperconvex sets.

\begin{proposition}
The intersection of a finite collection of sets in $\cH_k$ is 
again in $\cH_k$, but the intersection of an arbitrary 
collection of sets in $\cH_k$ is in $\cH_k$ if and only if the 
intersection is absorbing.
\end{proposition}
\begin{proof}
The result follows because the intersection of any collection 
of sets that satisfy the multi-variable von~Neumann inequality 
will again satisfy the multi-variable von~Neumann inequality.
\end{proof}

\begin{definition} Let $D \subseteq \ck$ be a subset of the 
closed unit polydisk,
then we let $\cH(D)$ denote the intersection of all elements in 
$\cH_k$ that contain $D.$
\end{definition}

\begin{proposition}
Let $D \subseteq \ck$ be a separating set contained in the 
closed unit 
polydisk. Then $\cH(D) \in \cH_k$ and if $H_1 \in \cH_k$ with 
$D \subseteq H_1$, 
then $\cH(D) \subseteq H_1.$
\end{proposition}

When $D$ is a subset of the closed unit polydisk that is 
separating then we refer to
$\cH(D)$ as the {\em IQ algebra ball generated by D.}
Again, we shall show that when $D$ is not separating, then 
$\cH(D)$ is not absorbing.

The following gives another way to realize $\cH(D).$

\begin{definition}
Let $D \subseteq \ck$ be a subset of the closed unit polydisk. 
For each $m$, 
for each polynomial $p$ in $m$ variables with $\|p\| \le 1$ and 
for each choice of
$m$ points, $w_1, \ldots, w_m$ in $D$, the vector $p(w_1, 
\ldots, w_m) \in \ck$ will lie
in the closed unit polydisk. We call the closure of the set of 
all such vectors in $\ck$,
the {\bf hyperconvex hull of D} and we denote it by $HC(D).$
\end{definition}

\begin{proposition}
Let $D \subseteq \ck$ be contained in the closed unit polydisk. 
Then $HC(D)$ is closed, absolutely
convex, $HC(D) \cdot HC(D) \subseteq HC(D)$ and is hyperconvex.
If $D$ is also separating, then $HC(D) \in \cH_k.$
\end{proposition}
\begin{proof}
An absolute convex combination of points in $HC(D)$ is the 
image of points in $D$ under the corresponding
absolute convex combination of the polynomials. Similarly, the 
product of two points in $HC(D)$ is
the image under the product of the corresponding polynomials. 
Finally, the fact that $HC(D)$ is hyperconvex
follows since the image of a collection of points in $HC(D)$ 
under a polynomial is the image of their
predecessors under a composition of polynomials.

It remains to show that if $D$ is separating, then $HC(D) \in 
\cH_k.$ 
It will be enough to show that $HC(D)$ is absorbing. Arguing as 
in the proof
of Proposition \ref{intersection}, we can show that a positive multiple of 
each basis vector is in $HC(D)$ and then
it follows that $HC(D)$ is absorbing.
\end{proof}

\begin{theorem}
Let $D \subseteq \ck$ be contained in the closed unit polydisk.
Then $HC(D) = \cH(D).$
\end{theorem}
\begin{proof}
If $D \subseteq H$ and $H \in \cH_k$, then necessarily $HC(D) 
\subseteq \cH(D).$
Thus, $HC(D) \subseteq H.$

If $D$ is separating then it follows that $HC(D) \in \cH_k$ and 
so $\cH(D) \subseteq HC(D).$
Now assume that $D$ is arbitrary. Without loss of generality, we 
may assume that $0 \in D$.
Let $\epsilon > 0$ be arbitrary and pick a vector $w=( w_1, 
\dots, w_k)$ with $w_i \ne w_j$ for $i \ne j$
and $w_i \le \epsilon$ for all $i.$
Let $D_1 = D \cup \{ w\}$, so that $D_1$ is separating.
Consequently, $\cH(D) \subseteq \cH(D_1) = HC(D_1).$

Finally, note that the only elements that are in $HC(D_1)$ but 
not in $HC(D)$ are polynomials that
involve the vector $w$ and other elements of $D$. Freezing all 
of the variables except the one involving
$w$ we see that we obtain a polynomial $p(z)$ in a single 
variable, with $\|p\| \le 1$.
Consequently, $|p(w_i) - p(0)| \le \epsilon$ and we deduce that 
every vector in $HC(D_1)$ is 
at most distance $k \epsilon$ from a vector in $HC(D)$.

Now assume that $v$ is not in $HC(D)$. Then for a small enough 
$\epsilon$, $v$ will not be in $HC(D_1)$.
But this latter set is the intersection of all balls in $\cH_k$ 
that contains $D_1$.
Hence, there will exist $H \in \cH_k$, such that $v$ is not in 
$H$, and $D \subseteq D_1 \subseteq H$.
Thus, $v$ is not in $\cH(D)$ and so $\cH(D) \subseteq HC(D).$
\end{proof}

\begin{corollary}
Let $D \subseteq \ck$ be contained in the closed unit polydisk. 
If $D$ is not separating,
then each of the sets $\cB(D), \cV(D), \cO(D)$ and $\cH(D)$ is 
not absorbing.
\end{corollary}
\begin{proof}
Since $\cB(D) \subseteq \cV(D) \subseteq \cO(D) \subseteq 
\cH(D)$, it will be enough to show that
$\cH(D)$ is not absorbing. But $\cH(D) = HC(D)$ and if there 
exists $i \ne j$ such that $w_i = w_j$
for every vector $w = (w_1, \ldots, w_k) \in D$, then it is 
easily seen that every vector in $HC(D)$
will be equal in the i-th and j-th components and so $HC(D)$ 
can not be absorbing.
\end{proof}

We close this section with another way to describe the sets 
$\cH(D)$ and $\cO(D)$ in the special case that $D= \{ v_1, 
\ldots, v_n \} \subseteq \ck$ is a separating set consisting of 
$n$ points. Recall that in this case setting $x_j = (w_{1,j}, 
\ldots, w_{n,j})$ where $v_i = ( w_{i,1}, \ldots, w_{i,k})$ 
defines $k$ distinct points in the closed polydisk in 
$\mathbb{C}^n$.

\begin{proposition} Let $D= \{ v_1, \ldots, v_n \} \subseteq 
\ck$ be a separating set consisting of $n$ points and let $x_1, 
\ldots, x_k$ be defined as above, then
$\cH(D) = \mathcal{D}( A(\mathbb{D}^n); x_1, \ldots, x_k).$
\end{proposition}
\begin{proof}
Given a polynomial $p$ in $n$ variables the i-th component of 
the vector $p(v_1, \ldots, v_n)$ is the number $p(x_i).$
Thus, $\mathcal{D}( A(\mathbb{D}^n); x_1, \ldots, x_k) 
\subseteq HC(D)= \cH(D).$ However, evaluating the $n$ 
coordinate functions at $x_1, \ldots, x_k$ shows that $D 
\subset \mathcal{D}( A(\mathbb{D}^n); x_1, \ldots, x_k)$ and 
since this latter set is hyperconvex, we have that $\cH(D) 
\subseteq   \mathcal{D}( A(\mathbb{D}^n); x_1, \ldots, x_k)$ 
and the result follows.
\end{proof}

A similar result holds for $\cO(D)$ if one first introduces the 
{\em universal operator algebra for n commuting contractions 
$A_u(\mathbb{D}^n)$ ,} as in \cite{Pau3}.
This algebra is the completion of the algebra of polynomials in 
$n$ variables where the norm of a polynomial is defined by 
taking the supremum of the norms of the operators defined by 
evaluating the polynomial at an arbitrary $n$ tuple of 
commuting operators on a Hilbert space.
By the fact that vonNeumann's inequality fails for 3 or more 
commuting contractions, we see that in general the norm of a 
polynomial in $A_u(\mathbb{D}^n)$ will generally be larger than 
its supremum norm over the polydisk.

Just as for a uniform algebra, given $x_1, \ldots, x_k$ in the 
closed polydisk, one may form the set
$$\mathcal{D}(A_u(\mathbb{D}^n); x_1, \ldots, x_k) = \{(f(x_1), 
\ldots, f(x_k)) : \|f\|_u \le 1 \}.$$
Applying either Cole's theorem \cite{BD} or the theorem of 
Blecher-Ruan-Sinclair, one finds that this set is the unit ball 
of an operator algebra norm on $\ck.$

\begin{lemma} Let $B \in \cO_k,$ let $v_1, \ldots, v_n \in B$ 
and let $f \in A_u(\mathbb{D}^n)$ with $\|f\|_u \le 1,$
then $f(v_1, \ldots, v_n) \in B.$
\end{lemma}
\begin{proof} It is enough to consider the case that $f$ is a 
polynomial, in which case the result follows from the 
factorization theory for universal operator algebra norms. See 
for example, \cite{Pau3} Corollary 18.2.
\end{proof}

\begin{proposition} Let $D= \{ v_1, \ldots, v_n \} \subseteq 
\ck$ be a separating set consisting of $n$ points and let $x_1, 
\ldots, x_k$ be defined as above, then
$\cO(D) = \mathcal{D}( A_u(\mathbb{D}^n); x_1, \ldots, x_k).$
\end{proposition}
\begin{proof}
Let $f \in A_u(\mathbb{D}^n)$ with $\|f\|_u \le 1.$ By the 
lemma, $(f(x_1), \ldots, f(x_k))= f(v_1, \ldots, v_n) \in 
\cO(D)$ and hence,
$\mathcal{D}( A_u(\mathbb{D}^n); x_1, \ldots, x_k) \subseteq 
\cO(D).$
 
Conversely, $\mathcal{D}( A_u(\mathbb{D}^n); x_1, \ldots, x_k) 
\in \cO_k$ and contains $D$ and hence contains $\cO(D).$
\end{proof}

\begin{corollary}
Let $D \subseteq \ck$ be contained in the closed unit polydisk. 
If $D$ consists of two or fewer points,
then $\cO(D) = \cH(D).$
\end{corollary}
\begin{proof}
By
Ando's theorem \cite{And}, we have that $A(\mathbb{D}^2) = 
A_u(\mathbb{D}^2)$ isometrically and hence $\cH(D) = 
\mathcal{D}(A(\mathbb{D}^2); x_1, x_2) = 
\mathcal{D}(A_u(\mathbb{D}^2); x_1, x_2)= \cO(D).$
\end{proof}

\section{Examples}

In this section we will present two examples. 
For the first example we will show that an example of Holbrook \cite{Hol} yields a $4$-idempotent
operator algebra acting on $\CC^4$ whose unit ball is not a hyperconvex set, i.e., an element of $\cO_4$ not in $\cH_4.$

Recall
that Cole and Wermer in \cite{CW2} show that $\ib$ is 
a semi-algebraic set when $A$ is the bidisk algebra. This leads naturally to the question of 
whether or not there exists a uniform algebra $A$ such that $\ib$ is 
not
a semi-algebraic set. For the second example we will construct
a $3$-idempotent operator algebra and conjecture that the unit ball  
determined by this algebra is not a semi-algebraic set. We have included 
a heuristic argument for why we believe the conjecture to be true.

\begin{example}
In this example we will use a result of J. Holbrook
to show that $\cH_4$ is a proper subset of $\cO_4$ and 
consequently
that $\cH_k \subsetneq \cO_k$ for all $k \geq 4$. 
\end{example}

Recall that in
1951 J. von Neumann proved that if $T$ is a contraction on a complex
Hilbert space and $p$ is polynomial in one variable (with complex
coefficients), then 
$$\| p(T) \| \leq \mathrm{sup} \{ |p(z)| :
|z| \leq 1 \}.
$$ 
T. Ando extended this result in 1963 by showing
that if $S$ and $T$ are commuting contractions and $p$ is 
a polynomial in two variables, then
$$
\| p(S,T) \| \leq \mathrm{sup} \{ |p(z,w)| : |z| \leq 1, \, |w| 
\leq 1 \}.
$$
In 1973, S. Kaijser and N. Th. Varopoulos explicitly describe 
$3$ commuting contractions $T_1$, $T_2$, and $T_3$ acting 
on a $5$-dimensional Hilbert space and a polynomial $p$ in 
three variables such that 
$$
\| p(T_1,T_2,T_3) \| > \|p \|_{\infty}.
$$
In 1991, Lotto and Steger proved a diagonalizable set of such 
contractions exist. Holbrook \cite{Hol} improved this result by lowering the dimension to 4.
\begin{theorem} \label{LottoSteger}
(Lotto-Steger and Holbrook) There are three commuting, diagonalizable 
contractions
$T_1$, $T_2$, and $T_3$ on $\CC^4$ and a polynomial $p$ in three
variables such that $\| p(T_1,T_2,T_3) \| > \|p\|_{\infty}$.
\end{theorem}

We use this fact 
to show that there exists a 
$4$-idempotent
operator algebra which is not an interpolation body (hyperconvex 
set).

Now let $T_1$, $T_2$, $T_3$, and $p$ be as in Theorem 
\ref{LottoSteger} acting on $\CC^4$
and choose an invertible $4 \times 4$ matrix $Q$ such that
$QT_jQ^{-1}$ for $j=1,2,3$ are the diagonal matrices
$$
QT_jQ^{-1} = \left(
\begin{array}{cccc}
w_1^j & 0 & 0 & 0\\
0 & w_2^j & 0 & 0 \\
0 & 0 & w_3^j & 0 \\
0 & 0 & 0 & w_4^j\\
\end{array} \right).
$$
Let $E_{ij}$ denote the canonical matrix
units and define idempotents 
$$E_i:= Q E_{ii} Q^{-1}$$
for
$i=1,2,3,4$ and let $\mathcal{A}:= \mathrm{span} \{E_1,E_2,E_3,E_4 \}$.
Next observe that $(w_1^j,w_2^j,w_3^j,w_4^j) \in 
\mathcal{D}(\mathcal{A})$
for $j=1,2,3$. However, 
$$p((w_1^1,...,w_4^1),(w_1^2,...,w_4^2),
(w_1^3,...,w_4^3)) \notin \mathcal{D}(\mathcal{A}).$$ Therefore,
$\mathcal{D}(\mathcal{A})$ fails to be hyperconvex and we have that 
$\cH_k \subsetneq \cO_k$ for all $k \geq 4$.


\begin{example} In this example we will inductively
construct a Schur ideal (non-trivial and bounded) and we conjecture that 
the ``perp'' of 
this Schur ideal is not a semi-algebraic subset of 
$\mathbb{R}^3$.
Hence, if the conjecture is true, then by Proposition \ref{schuralgebra} we will have that  
there is a $3$-idempotent operator algebra such that the
ball of this operator algebra is not a semi-algebraic set.
\end{example}

First we will inductively construct the appropriate Schur ideal.

\begin{lemma} Let $a,c >0$. Then for $x,y \in \mathbb{R}$ we 
have 
that $(0,x,y) \in \{P_{a,c} \}^{\p}$ if and only if
$$
y^2 \leq \frac{ac-(ac+c)x^2}{(ac+a)-(ac+a+c)x^2},
$$
where $P_{a,c}$ is the following $3\times 3$ positive definite 
matrix:
$$
\left(
\begin{array}{ccc}
1 & 1 & 1 \\
1 & a+1 & 1 \\
1 & 1 & c+1 \\
\end{array}
\right).
$$
\end{lemma}

\begin{proof}
By definition $(0,x,y) \in \{ P_{a,c} \}^{\p}$ if and only if 
the following matrix
\begin{eqnarray}
\label{matrix1}
\left(
\begin{array}{ccc}
1 & 1 & 1\\
1 & 1-x^2 & 1-xy\\
1 & 1-xy & 1 -y^2\\
\end{array}
\right) \ast P_{a,c}
\end{eqnarray}
is positive semi-definite. 
The matrix in (\ref{matrix1}) is positive semi-definite if and 
only if 
\begin{eqnarray}
(a-(a+1)x^2)(c-(c+1)y^2)- x^2y^2 \geq 0  \text{ and } 
\label{eq1}\\
x^2 \leq \frac{a}{a+1} \text{ and } y^2 \leq \frac{c}{c+1} 
\label{eq2}.
\end{eqnarray}
by applying the Cholesky algorithm, \cite{HJ}.
Using (\ref{eq2}) we have that $ac-(ac+c)x^2 \geq 0$ and 
$ac+a-(ac+a+c)x^2 \geq 0$. Therefore we have that the matrix in 
(\ref{matrix1})
is postive semi-definite if and only if 
$$
y^2 \leq \frac{ac-(ac+c)x^2}{(ac+a)-(ac+a+c)x^2}.
$$
\end{proof}

Let $u=x^2, v=y^2$ and let $f_{a,c}(u):=\frac{ac-(ac+c)u}{(ac+a)-(ac+a+c)u}$. Note
that $f_{a,c}$ is decreasing
with $f_{a,c}(0)=\frac{c}{c+1}$ and $f_{a,c}\left( \frac{a}{a+1} 
\right)=0$.

Next observe that if $a_1 < a$ and $c_1 > c$, then
$\frac{a_1}{a_1+1} < \frac{a}{a+1}$ and
$\frac{c_1}{c_1+1} > \frac{c}{c+1}$. Thus the graphs 
of $f_{a,c}$ and $f_{a_1,c_1}$ intersect at some point
$\mu_1< \frac{a_1}{a_1 +1}$.

For each function in the family, $f^{\prime}_{a,c}(u) = \frac{-ac}{(ac-(ac+c)u)^2}$ is also a 
monotone decreasing function of $u.$
Since, $\displaystyle{ \lim_{c_1 \to +\infty} f^{\prime}_{a_1,c_1}(0) = -\infty},$ by chosing $c_1$ 
sufficiently large we can guarantee that at $\mu_1$, we have
$f^{\prime}_{a_1,c_1}(\mu_1) < f^{\prime}_{a_1,c_1}(0) < f^{\prime}_{a,c}(\mu_1).$

Now choose any $a_2 < a_1$ such that $\mu_1 < \frac{a_2}{a_2 +1}$ and then choose $c_2 > c_1$ 
such that the point of intersection, $\mu_2,$ of $f_{a_1,c_1}$ and $f_{a_2,c_2}$ satisfies $\mu_1 < \mu_2.$
To see that this can be done, note that for fixed $a_2,$ as $c_2$ increases the point of intersection of 
the two curves moves to the right continuously and approaches $\frac{a_2}{a_2+1}$ in the limit. 
This also allows us to choose $c_2$ so that $|\frac{a_2}{a_2+1} - \mu_2| \le \frac{1}{2}|\frac{a_1}{a_1+1} - \mu_1|.$

Now if we 
look at the set
$$
\{ (u,v): v \leq f_{a,c}(u), \, v \leq f_{a_1,c_1}(u), \, 
v \leq f_{a_2,c_2}(u) \}
$$ 
it will have two non-differentiable corners. 

Now inductively choose $\{ a_n\}$ and $\{ c_n \}$ such that
the $a_n \searrow$ and $c_n \nearrow$ and so that if 
$\mu_n$ is the point of intersection of $f_{a_n,c_n}$
with $f_{a_{n+1},c_{n+1}}$, then we have that $\mu_n \nearrow,$ with $|\frac{a_{n+1}}{a_{n+1}+1} - \mu_{n+1}| \le \frac{1}{2}|\frac{a_n}{a_n +1} - \mu_n|$
and $\mu_n < \frac{a_m}{a_m +1}$ for all $m,n \in \mathbb{N}$. 
Note that $\displaystyle{ \lim_{n \to \infty} \mu_n = \lim_{n \to \infty} \frac{a_n}{a_n + 1}}$ and we call this common limit $\mu.$
Set $a_0 =a, b_0=b$ and $\mu_0=0.$
 
\begin{conjecture}
Let $a_n, \, c_n >0$ be chosen as above. Then
 $\{ P_{a_n,c_n} : n \ge 0 \}^{\p}$ is a non-algebraic subset of $\cO_4.$
\end{conjecture}

Our only obstruction to proving the above conjecture is a problem concerning semi-algebraic sets that seems likely to be true, but which we have been unable to prove or find in the literature.

\begin{conjecture}
Let $f: [0,1] \to \mathbb{R}$ be continuous, non-negative and set $C= \{ (x,y): 0 \le x \le 1, 0 \le y \le f(x) \}.$ If $f$ is non-differentiable at infinitely many points, then $C$ is not semi-algebraic.
\end{conjecture}

\begin{proposition} If the above conjecture concerning semi-algebraic sets is true, then there exists $B \in \cO_3$ that is not semi-algebraic.
\end{proposition}

\begin{proof} We let $B= \{ P_{a_n,c_n} : n \ge 0 \}^{\p}$ denote the set in the earlier conjecture.

By our earlier results, it is clear that $B$ is in $\cO_4.$
So it remains to show that this set is not semi-algebraic.

Now for $x$ and $y$ real, $(0,x,y) \in B$,
if and only if $y^2 \leq f(x^2)$ where
$$
f(u) := \inf \{ f_{a_n,c_n}(u) : n \in \mathbb{N} \}.
$$
Also, for $\mu_n \le u \le \mu_{n+1},$ we have that $f(u) = f_{a_n,c_n}(u).$

Now if $B$ was semi-algebraic, then $C= \{ (x,y) \in \mathbb{R}^2 : (0,x,y) \in B \}$ is semi-algebraic, by 
the Tarski-Seidenberg theorem\cite{CW2}. Loosely speaking, Tarski-Seidenberg says that if a subset $X \subseteq \mathbb{R}^{n+1}$ is a semi-algebraic set, then 
any coordinate projection of $X$ onto $\mathbb{R}^n$ is also a semi-algebraic set.

But by our construction, $C$ is the region under the graph of a function $f$ that has a countable collection of points of non-differentiability.
By the second conjecture, 
this set can not be semi-algebraic.

\end{proof}

Cole and Wermer\cite{CW2} has an appendix which includes many of the important theorems on semi-algebraic sets and serves as a nice introduction to this area.

\end{document}